\documentclass[a4paper, 11pt, reqno]{amsart}
\usepackage{amssymb, amsfonts}
\usepackage{bookmark, bm, cleveref, enumitem, hyperref}
\crefname{enumi}{equality}{}
\crefname{equation}{}{}
\creflabelformat{equation}{(#2#1#3)}

\usepackage{tikz}
\usetikzlibrary{positioning, arrows, petri}

\numberwithin{equation}{section}

\hypersetup{pdfstartview={FitH}}

\newtheorem{tm}{Theorem}
\newtheorem{lm}[tm]{Lemma}
\newtheorem{prop}[tm]{Proposition}
\newtheorem{kor}[tm]{Corollary}
\newtheorem{remark}[tm]{Remark}
\newenvironment{dokaz}{\noindent\emph{Proof:}\ }{\hfill $\blacksquare$}

\newcommand{\Z}{{\mathbb Z}}
\newcommand{\N}{{\mathbb N}}
\newcommand{\C}{{\mathbb C}}

\newcommand{\g}{{\mathfrak g}}
\newcommand{\gk}{\hat{{\mathfrak g}}}
\newcommand{\gt}{\tilde{{\mathfrak g}}}
\newcommand{\nt}{\tilde{{\mathfrak n}}}
\newcommand{\h}{{\mathfrak h}}
\newcommand{\n}{{\mathfrak n}}

\newcommand{\Gamt}{\tilde{{\Gamma}}}
\newcommand{\mm}{{\bf m}}

\newcommand{\gdva}{{\gamma_2}}
\newcommand{\gtri}{{\gamma_3}}
\newcommand{\gcet}{{\gamma_4}}
\newcommand{\gll}{{\gamma_\ell}}

\newcommand{\gllu}{{\gamma_\llu}}

\newcommand{\gdu}{{\gamma_\dvau}}
\newcommand{\gtu}{{\gamma_\triu}}
\newcommand{\gcu}{{\gamma_\cetu}}

\newcommand{\xdu}{{x_\dvau}}
\newcommand{\xtu}{{x_\triu}}
\newcommand{\xcu}{{x_\cetu}}
\newcommand{\llu}{{\underline{\ell}}}
\newcommand{\guu}{{\underline{\gamma}}}
\newcommand{\dvau}{{\underline{2}}}
\newcommand{\triu}{{\underline{3}}}
\newcommand{\cetu}{{\underline{4}}}

\begin{document}

\title{Generating sets of standard modules for $D_4^{(1)}$}
\author{Ivana~Baranovi\'c}
\thanks{\textsc{Ivana~Baranovi\'c}, \texttt{baranovic.ivana@yahoo.com}\vspace{1mm}}
\author{Miroslav~Jerkovi\'c}
\thanks{\textsc{Miroslav Jerkovi\'c}, \texttt{mjerkov@fkit.unizg.hr} \\ \indent{\textsc{University of Zagreb Faculty of Chemical Engineering and Technology, Zagreb, Croatia}}\vspace{1mm}}
\author{Goran~Trup\v cevi\'c}
\thanks{\textsc{Goran Trup\v cevi\'c}, \texttt{goran.trupcevic@ufzg.hr} \\ \indent{\textsc{University of Zagreb Faculty of Teacher Education, Zagreb, Croatia}}}
\begin{abstract}
	Let $\tilde\g$ be an affine Lie algebra of type $D_4^{(1)}$ and $L(\Lambda)$ its standard module of level $k$ with highest weight vector $v_{\Lambda}$. We define Feigin-Stoyanovsky's type subspace as
	\[
		W(\Lambda)=U(\tilde\g_{1})\,v_{\Lambda},
	\]
	where $\tilde\g=\tilde\g_{-1}\oplus\tilde\g_{0}\oplus\tilde\g_{1}$ is a $\Z$-gradation of $\tilde\g$ associated with a $\Z$-gradation $\g=\g_{-1}\oplus\g_{0}\oplus\g_{1}$. Using vertex operator relations, we reduce the Poincaré-Birkhoff-Witt spanning set of $W(\Lambda)$, and describe it in terms of difference and initial conditions. The spanning set of the whole standard module $L(\Lambda)$ can be obtained as a limit of the spanning set for $W(\Lambda)$.
\end{abstract}

\maketitle

\section{Introduction}

Let $\g$ be a finite-dimensional Lie algebra and $\tilde\g$ its corresponding affine Lie algebra
\[
	\tilde\g=\g\otimes\C[t,t^{-1}]\oplus\C c \oplus\C d,
\]
where $c$ is the canonical central element and $d$ the degree operator (cf. \cite{K}). Affine Lie algebras have a vertex operator construction and this construction has been used by J.~Lepowsky and R.~Wilson \cite{LW1,LW2} to give a Lie-theoretic interpretation of Rogers-Ramanujan type combinatorial identities. This approach was then continued in the works of J.~Lepowsky, M.~Primc, A.~Meurman, and others (cf. \cite{LP,MP}). An important part of the program was to find monomial bases of standard modules (or some of their subspaces) for affine Lie algebras.

In \cite{FS}, B.~Feigin and A.~Stoyanovsky introduced and studied a principal subspace for $\widetilde{\mathfrak{sl}}(\ell+1,\C)$. This subspace is defined as $W(\Lambda)=U(\tilde\n)v_\Lambda$, where $\n$ is the nilradical of a Borel subalgebra of $\tilde{\mathfrak{sl}}(\ell+1)$:
\[
	\tilde\n=\n\otimes\C[t,t^{-1}],
\] and $\Lambda$ and $v_\Lambda$ are the highest weight and a highest weight vector of $L(\Lambda)$, respectively. This definition can be extended to an arbitrary highest weight module for an affine Lie algebra.

Galin Georgiev (cf. \cite{Ge}) studied further principal subspaces and constructed a certain class of combinatorial bases for representations of $\widetilde{\mathfrak{sl}}(\ell+1,\C)$. He used intertwining operators from \cite{DL} to prove linear independence of constructed bases. These operators have also been used by S.~Capparelli, J.~Lepowsky, and A.~Milas (see \cite{CLM1,CLM2}) to obtain the Rogers-Ramanujan and Rogers-Selberg recursions for characters of principal subspaces for $\widetilde{\mathfrak{sl}}(\ell+1,\C)$. Along with C.~Calinescu (see \cite{C}), the continuation of this work led to Rogers-Ramanujan-type recursions for the principal subspaces of level $1$ standard representations of $\widetilde{\mathfrak{sl}}(\ell+1,\C)$ and for certain principal subspaces of the higher-level standard modules for $\widetilde{\mathfrak{sl}}(3,\C)$. Solving these recursions gave character formulas of the principal subspaces.

Spaces similar to the principal subspaces have been studied by M.~Primc in \cite{P1,P2,P3}, and bases of these subspaces were then used to construct bases of standard modules. In \cite{P3}, M.~Primc used the Capparelli-Lepowsky-Milas approach to provide a simpler proof for the linear independence of bases constructed in \cite{P1}. This approach has been shown to be fruitful for other Feigin-Stoyanovsky's subspaces for $\widetilde{\mathfrak{sl}}(\ell+1,\C)$ as well; in \cite{T1}, G.~Trupčević found the combinatorial bases of $W(\Lambda)$ for all level one standard modules for $\widetilde{\mathfrak{sl}}(\ell+1,\C)$. He also extended the proof to all higher level modules in \cite{T2}. In \cite{J}, M.~Jerković obtained the exact sequences of Feigin-Stoyanovsky's type subspaces at fixed level $k$ using the known description of combinatorial bases for Feigin-Stoyanovsky's type subspaces of standard modules for the affine Lie algebra $\widetilde{\mathfrak{sl}}(\ell+1,\C)$. This led to systems of recurrence relations for formal characters of those subspaces.

In this note we use the same approach to continue the work started in \cite{Ba}. Using the vertex operator relations we find the spanning set for Feigin-Stoyanovsky's type subspaces of standard modules of level $k$ for algebra $D_{4}^{(1)}$. We then show that the spanning set of the whole standard module $L(\Lambda)$ can be obtained as a limit of a spanning set for $W(\Lambda)$.

The paper is organized as follows: In Sections \ref{S: Affine} and \ref{S: FS} we give the setting and introduce the notion of Feigin-Stoyanovsky's subspace, together with a linear ordering on the set of monomials that span this subspace. In Section~\ref{S: VOA} we review the well-known Frenkel-Kac-Segal construction of the basic $\tilde\g$-modules, define the simple current operator and review vertex operator relations. Next, in Section~\ref{S: ICDC} we show that monomials satisfying the difference and initial conditions form a spanning set for $W(\Lambda)$. Finally, in Section~\ref{S: Spanning} we use this result to obtain a generating set for the whole standard module $L(\Lambda)$.

\medskip
\noindent\textit{2000 Mathematics Subject Classification.} Primary 17B67; Secondary 17B69, 05A19.

\noindent\textit{Acknowledgement.} We sincerely thank Mirko Primc for his help on the topic and valuable suggestions.

\section{Affine Lie algebras and basic modules} \label{S: Affine}

For $\ell\geq 4$ let ${\mathfrak g}$ be a complex simple Lie algebra of type $D_\ell$ and $\h$ its Cartan
subalgebra. The corresponding root system $R$ may be realized in $\mathbb R\sp\ell$. Fix the basis $\Pi=\{\alpha_1,\dots,\alpha_\ell\}$ of $R$ where 
$$\alpha_1=\epsilon_1-\epsilon_2, \dots, \alpha_{\ell-1}=\epsilon_{\ell-1}-\epsilon_\ell, \alpha_\ell=\epsilon_{\ell-1}+\epsilon_\ell,$$
where $\epsilon_1,\dots,\epsilon_\ell$ is the canonical basis of $\mathbb R\sp\ell$.
Then we have a root space decomposition ${\mathfrak g}={\mathfrak h}+\sum_{\alpha\in R}\mathfrak g_\alpha$ and a triangular decomposition $ \g=\n_-+ \h + \n_+$ of $\g$. 
Let $\theta=\epsilon_1+\epsilon_2$ be the maximal root. 
Denote by $\langle\,,\,\rangle$ the normalized Killing form such that $\langle\theta,\theta\rangle=2$. We identify $\mathfrak h$ and $\mathfrak h\sp*$ via the Killing form. 
For each root $\alpha$ we fix a root vector $x_\alpha\in\g_\alpha$.

Denote by $\omega_1,\dots,\omega_\ell$ the fundamental weights and define, for convenience, $\omega_0=0$. We fix, as in \cite{Ba}, weight $\omega=\omega_1=\epsilon_1$. Note that $\omega$ is a minuscule weight, i.e.
$$\langle\omega,\alpha\rangle\in\{-1,0,1\}\quad\text{for all}\quad \alpha\in R.
$$
This gives the induced $\mathbb Z$-gradation of $\g$
$$
	\g =\g_{-1} + \g_0 + \g_1, \qquad
	\g_0 = \h +
	\sum_{\langle\omega,\alpha\rangle=0}\, \g_\alpha,
	\qquad
	\g_{\pm1} = 
	\sum_{\langle\omega,\alpha\rangle= \pm 1}\, \g_\alpha.
$$
The subalgebras $\g_{1}$ and $\g_{-1}$ are commutative subalgebras and $\g_0$ acts on them by adjoint action. The subalgebra $\g_0$ is reductive with semisimple part $\mathfrak{l}_0=[\g_0,\g_0]$ being a simple algebra of type $A_3$ if $\ell=4$ and of type $D_{\ell-1}$ if $\ell\geq 5$.

Set $\Gamma = \{\,\alpha \in R \mid \langle\omega,\alpha\rangle = 1\}$. Then 
$$\Gamma = \{\gdva,\dots,\gll,\gllu, \dots, \gdu\},$$
where $\gamma_i =\epsilon_1+\epsilon_i$ and 
$\gamma_{\underline{i}} =\epsilon_1-\epsilon_i$, 
for $i=2,\dots,\ell$.
We say that $\Gamma$ is {\it the set of colors} and we write
$ x_i =x_{\gamma_i}$, 
$x_{\underline{i}} = x_{\gamma_{\underline{i}}}$. 
Notice that $\gdva$ is the maximal root $\theta$.

Denote by $\gt$ the associated affine Lie algebra of type $D_\ell\sp{(1)}$,
$$ \tilde{\mathfrak g} = \mathfrak{g} \otimes
	\mathbb{C}[t,t^{-1}] + \mathbb{C}c + \mathbb{C} d,
$$
where $c$ is the canonical central element and $d$ the degree operator (cf. \cite{K}). Set
$
	x(j)=x\otimes t^{j}
$
for $x\in{\mathfrak g}$ and $j\in\mathbb Z$ and denote by $ x(z)=\sum_{n\in\mathbb Z}
	x(n) z^{-n-1}$ a formal Laurent series in formal variable $z$. 

Set
$\h^e=\h+\C c+\C d$,
$\tilde{\n}_\pm=\g\otimes t^{\pm 1}\C [t^{\pm 1}] + \n_\pm,$

We have the triangular decomposition: $
	\gt=\tilde{\n}_-+ \h^e + \tilde{\n}_+$, and the induced $\Z$-gradation with respect to $\omega$: 
$$
	\gt=\gt_{- 1}+\gt_{0}+\gt_{ 1},
$$
$
	\gt_0 = {\mathfrak g}_0\otimes\C [t,t^{-1}] + \C c + \C d$, $\gt_{\pm 1} = {\mathfrak g}_{\pm 1}\otimes\C [t,t^{-1}]$.
The subspace $\gt_1\subset \gt$ is a commutative subalgebra and $\g_0$ acts on
$\gt_1$ by adjoint action.

Let $\tilde{\Pi}=\{\alpha_0,\alpha_1,\dots,\alpha_\ell\}\subset(\h^e)^*$ be the set of simple roots of $\gt$.
Denote by $\Lambda_0,\dots,\Lambda_\ell\in(\h^e)^*$ the fundamental weights of $\gt$ defined by $\langle \Lambda_i,\alpha_j\rangle=\delta_{ij}$
and $\Lambda_i(d)=0$, for $i=1,\dots,\ell$ where $\langle\cdot,\cdot\rangle$ is the standard extension of Killing form onto $\h^e$.

Denote by $L(\Lambda)$ a standard $\tilde{\mathfrak g}$-module, i.e. irreducible highest weight module with dominant integral highest weight
$$
	\Lambda=k_0 \Lambda_0+k_1 \Lambda_1+\dots+k_\ell \Lambda_\ell,
$$
$k_i\in\Z_+$ for $i=0,\dots,\ell$, 
and by $v_\Lambda$ a fixed highest weight vector of $L(\Lambda)$. The central element $c$ acts on $L(\Lambda)$ as a scalar 
$$
	k=\Lambda(c)=k_0 +k_1 +2k_2+\dots+2k_{\ell-2}+k_{\ell-1}+k_\ell,
$$
which is called the level of the module $L(\Lambda)$.
Note that modules $L(\Lambda_0)$, $L(\Lambda_1)$, $L(\Lambda_{\ell-1})$ and $L(\Lambda_\ell)$ are of level $1$ whilst $L(\Lambda_2),\dots,L(\Lambda_{\ell-2})$
are level $2$ modules.

\section{Feigin-Stoyanovsky's type subspaces} \label{S: FS}

For each integral dominant $\Lambda$ define a {\it Feigin-Stoyanovsky's type subspace}
$$
	W(\Lambda)=U(\gt_1)\cdot v_\Lambda\subset L(\Lambda).
$$
If we set $\gt_1^-=\gt_1\cap \nt_-$, then
$$
	W(\Lambda)=U(\gt_1^-)\cdot v_\Lambda.
$$
By Poincare-Birkhoff-Witt theorem, there is a spanning set of $W(\Lambda)$ consisting of {\it monomial vectors}
\begin{equation} \label{E: PBW base def}
	\{x_{\gamma_1}(-n_1)
	x_{\gamma_2}(-n_2)\cdots x_{\gamma_t}(-n_t) v_\Lambda\,|\,
	t\in\Z_+,\gamma_i\in \Gamma, n_i\in\N\}.
\end{equation}
Elements of this spanning set can be identified with monomials from $U(\gt_1^-)=S(\gt_1^-)$. We will sometimes refer to the elements from $\gt_1^-$
as {\it variables}, {\it elements} or {\it factors} of a monomial.

Introduce a linear order on monomials in the following way.
First, define a linear order on the set of colors $\Gamma$:
$$\gdva > \gtri > \dots > \gll > \gllu > \dots > \gdu.$$
On the set of variables $\Gamt=\{x_\gamma(n) \,|\,\gamma\in\Gamma, n\in \Z \}$ define a linear order by
$x_\alpha(i)<x_\beta(j)$ if $i<j$ or
$i=j,\, \alpha<\beta$. 
Finally, to define order on monomials from $S(\gt_1^-)$ assume that variables are sorted descendingly from right
to left. The order $<$ on the set of monomials is
defined as a lexicographic order, where we compare variables from 
right to left (from the greatest to the lowest one).

Order $<$ is compatible with multiplication (see \cite{P1, T1}):
\begin{equation} \label{OrdMult_rel}
	\textrm{if} \quad {\bf m} < {\bf m}'\quad \textrm{then}\quad {\bf m} {\bf m}_1 <{\bf m}' {\bf m}_1\quad 
\end{equation}
for monomials ${\bf m}, {\bf m}', {\bf m}_1\in\C[\Gamt]$.

We also define the {\em length} of a monomial to be the number of its factors.

\section{Vertex operator construction} \label{S: VOA}

For the rest of the paper, we restrict our considerations to the case of $\mathfrak{g}$ of the type $D_4$. In this case, the set of colors is $\Gamma=\{ \gdva,\gtri,\gcet,\gcu,\gtu,\gdu \}$. The fundamental weights for $\gt$ are $\Lambda_0,\dots,\Lambda_4$, and the corresponding standard modules $L(\Lambda_i)$ are of level $1$ for $i=0,1,3,4$ and of level $2$ for $i=2$.

By Frenkel-Kac-Segal construction (see \cite{LL} for details), standard modules can be realized on the tensor product $V_P=M(1)\otimes\C [P]$ of the Fock space $M(1)$ for the homogeneous Heisenberg subalgebra and the group algebra $\C [P]$ of the weight lattice with the basis $e^\lambda, \lambda\in P$. The action of $\gt$ is given by the vertex operator formula:
$$
	x_\alpha(z)=E^-(-\alpha,z)E^+(-\alpha,z)e^\alpha \epsilon(\alpha,\cdot) z^\alpha,
$$
where $z^\lambda \cdot e^\mu=e^\mu z^{\langle \lambda,\mu \rangle}$, $e_\alpha=e^\alpha \epsilon(\alpha,\cdot)$ is a product of the multiplication operator $e^\alpha$ and a $2$-cocycle $\epsilon(\cdot,\cdot)$, and 
$ E^{\pm}(\alpha,z)=\exp \left(\sum_{m\geq 1}\alpha(\pm m) {z^{\mp m}}/{\pm m}\right)\otimes 1$.

The basic module $L(\Lambda_0)=V_Q=M(1)\otimes \C [Q]$ is a vertex operator algebra, and $L(\Lambda_i)=M(1)\otimes e^{\omega_i}\C [Q]$ for $i=1,3,4$ are modules for this algebra. Level $2$ standard module $L(\Lambda_{2})$ can be realized as a submodule of the tensor product of two copies of a spinor representation, 
namely $L(\Lambda_{3}) \otimes L(\Lambda_{3})$ (and similarly for $L(\Lambda_{4}) \otimes L(\Lambda_{4})$) (see \cite{Ba}).
Standard modules of level $k>1$ can be realized as submodules of the tensor product 
of level $1$ and level $2$ standard modules:
$$L(\Lambda)\subset L(\Lambda_{0})^{k_0} \otimes L(\Lambda_{1})^{k_1} \otimes L(\Lambda_{2})^{k_2} \otimes L(\Lambda_{3})^{k_3} \otimes L(\Lambda_{4})^{k_4} \subset V_P^{\otimes k}.
$$
By taking tensor products of vertex operators $V_P^{\otimes k}$ becomes a module for the vertex operator algebra $V_Q^{\otimes k}$ (see \cite{DL}).

We now recall some relations between vertex operators from \cite{Ba} and extend them to the level $k$ case.

\begin{lm}\label{DC1_lm}
	For $\gamma,\delta,\tau\in \Gamma$
	\begin{align}
		x_\delta(-1) x_\gamma(-1) v_{\Lambda_{0}}            & = 0 \qquad \textrm{if}\qquad \delta\neq\guu, \label{voa_rel1}                                                \\
		x_\guu(-1) x_\gamma(-1) v_{\Lambda_{0}}              & = C x_{\underline{\delta}}(-1) x_\delta(-1) v_{\Lambda_{0}} = C' e^{2\omega} v_{\Lambda_0}, \label{voa_rel2} \\
		x_\tau(-1) x_\delta(-1) x_\gamma(-1) v_{\Lambda_{0}} & = 0, \label{voa_rel3}
	\end{align}
	where $C,C'\in\C$ are some nonzero constants.
\end{lm}

This gives relations between vertex operators for level $1$ modules (cf. \cite{DL, LL}):

\begin{tm} \label{voarel1_tm} Let $\gamma,\delta\in\Gamma$. Then on $L(\Lambda_i)$, $i=0,1,\ell-1,\ell$, the following holds
	\begin{eqnarray*}
		& & \vspace{-2ex} x_\delta(z) x_\gamma(z) = 0 \qquad \textrm{if}\qquad \delta\neq\guu, \\
		& & \vspace{-2ex} x_\guu(z) x_\gamma(z) = C x_{\underline{\delta}}(z) x_\delta(z),  \\
		& & \vspace{-2ex} x_\delta(z) x_\guu(z) x_\gamma(z) = 0, 
	\end{eqnarray*}
	for some nonzero constant $C\in\C$.
\end{tm}

The same principle can be applied for higher-level modules. From relations \cref{voa_rel1,voa_rel2,voa_rel3} it follows that
\begin{equation*}
	x_{\delta_1}(-1) x_{\delta_2}(-1) \cdots x_{\delta_{k+1}}(-1) v_{\Lambda_0}^{\otimes k} =0,
\end{equation*}
if the set $\{\delta_1,\delta_2,\dots,\delta_{k+1}\}$ does not contain any pair $\gamma, \guu$. Again, this gives relations between
vertex operators for level $k$ modules $L(\Lambda)$:
\begin{equation} \label{VOAdc_rel}
	x_{\delta_1}(z) x_{\delta_2}(z) \cdots x_{\delta_{k+1}}(z) =0.
\end{equation}

We also recall relations that gave rise to the initial conditions for level $1$ and $2$ modules from \cite{Ba}.

\begin{lm} \label{IC1_lm}
	For level $1$ modules $L(\Lambda_1)$, $L(\Lambda_3)$ and $L(\Lambda_4)$ the following relations hold:
	\begin{eqnarray*} 
		& & \vspace{-2ex} x_\gamma(-1) v_{\Lambda_{1}} = 0 \qquad \textrm{for}\qquad \gamma\in\Gamma,\\
		& & \vspace{-2ex} x_\gamma(-1) v_{\Lambda_{3}} = 0 \qquad \textrm{for}\qquad \gamma\in\{2,3,\cetu\},\\
		& & \vspace{-2ex} x_\gamma(-1) v_{\Lambda_{4}} = 0 \qquad \textrm{for}\qquad \gamma\in\{2,3,4\}.
	\end{eqnarray*}
	For level $2$ module $L(\Lambda_2)$	the following relations hold		
	\begin{eqnarray*}
		& & \vspace{-2ex} x_2(-1) v_{\Lambda_{2}} = 0,\\
		& & \vspace{-2ex} x_\tau(-1) x_\delta(-1) v_{\Lambda_{2}} = 0 
	\end{eqnarray*}
	for $ \tau \leq \delta$ and $\tau\neq\underline{\delta}$.
\end{lm}

Define a simple current operator $[\omega] : V_P\to V_P$ by $[\omega] =e^\omega \epsilon(\cdot,\omega)$ whose restrictions 
\begin{eqnarray*}
	& & \hspace{-2ex} L(\Lambda_0) \xrightarrow{[\omega]} L(\Lambda_{1}) \xrightarrow{[\omega]}
	L(\Lambda_0)\\
	& & \hspace{-2ex} L(\Lambda_3) \xrightarrow{[\omega]} L(\Lambda_{4}) \xrightarrow{[\omega]}
	L(\Lambda_3)
\end{eqnarray*}
are simple current operators on level $1$ standard modules. The following relations hold:
\begin{align}
	[\omega] v_0 & = C_1 v_{1} \label{SCrel1}                   \\
	[\omega] v_1 & = C_2 x_2(-1) x_\dvau(-1) v_0 \label{SCrel2} \\
	[\omega] v_3 & = C_3 x_\cetu(-1) v_4 \label{SCrel3}         \\
	[\omega] v_4 & = C_4 x_4(-1) v_3 \label{SCrel4}
\end{align}
for some nonzero constants $C_1,C_2,C_3,C_4$.

The following commutation relation holds 
\begin{equation*} 
	x_\alpha(n)[\omega]=[\omega] x_\alpha (n+\alpha(\omega)), \qquad \alpha\in R, \quad n\in\mathbb Z,
\end{equation*}
(cf. \cite{DLM} and \cite{Ga}). 
For $\gamma \in \Gamma$ one has
$$	x_\gamma(-n-1) [\omega] = [\omega] x_\gamma(-n).$$
For a monomial ${\bf m}\in U(\gt_1)$, denote by ${\bf m}^\pm$ a monomial obtained by raising/decreasing degrees in ${\bf m}$ by $1$. By
${\bf m}^{\pm m}$ denote a monomial obtained by raising/decreasing degrees in ${\bf m}$ by $m$. Then
$$
	{\bf m} [\omega] = [\omega]{\bf m}^+. 
$$

In the level $k>1$ case, 
use tensor products of level $1$ simple current operators
$$[\omega]=[\omega]\otimes\dots\otimes [\omega].$$
In particular, $L(\Lambda_2) \xrightarrow{[\omega]} L(\Lambda_{2})$. The following relation holds
\begin{equation} \label{SCrel5}
	[\omega] v_2 = C\cdot x_3(-1) x_\triu (-1) v_2,
\end{equation}
for some nonzero constant $C$.

\section{Difference and initial conditions} \label{S: ICDC}

For a monomial $\mm$ and a given $j\in\N$ let $a_t$ (resp. $a_{\underline{t}}$) denote the number of $x_t(-j)$ (resp. $x_{\underline{t}}(-j)$) factors in $\mm$. 
Furthermore, let $b_t$ (resp. $b_{\underline{t}}$) denote the number of $x_t(-j-1)$ (resp. $x_{\underline{t}}(-j-1)$) factors in $\mm$. We say that the monomial $\mm$
satisfies \emph{difference conditions} for level $k$ if the following frequency conditions hold:
\begin{enumerate}[label=(\roman*)]
	\item \label{DC1} $b_2+ a_\dvau+ a_\triu+ a_\cetu+ a_4+ a_3 \leq k$,
	\item \label{DC2} $b_2+ a_\triu+ a_\cetu+ a_4+ a_3+ a_2 \leq k$,
	\item \label{DC3} $b_3+ b_2+ a_\dvau+ a_\triu+ a_\cetu+ a_4 \leq k$,
	\item \label{DC4} $b_3+ b_2+ a_\dvau+ a_\cetu+ a_4+ a_3 \leq k$,
	\item \label{DC5} $b_4+ b_3+ b_2+ a_\dvau+ a_\triu+ a_4 \leq k$,
	\item \label{DC6} $b_\cetu+ b_3+ b_2+ a_\dvau+ a_\triu+ a_\cetu \leq k$,
	\item \label{DC7} $b_\cetu+ b_4+ b_3+ b_2+ a_\dvau+ a_\triu \leq k$,
	\item \label{DC8} $b_\triu+ b_\cetu+ b_4+ b_2+ a_\dvau+ a_\triu \leq k$,
	\item \label{DC9} $b_\triu+ b_\cetu+ b_4+ b_3+ b_2+ a_\dvau \leq k$,
	\item \label{DC10} $b_\dvau+ b_\triu+ b_\cetu+ b_4+ b_3+ a_\dvau \leq k$,
\end{enumerate}
for every $j\in\N$.

We recall results on Feigin-Stoyanovsky's type subspace of level $1$ and level $2$ standard modules for algebra $D_4^{(1)}$ from \cite{Ba}.
Relations between vertex operators from Theorem \ref{voarel1_tm} enabled us to reduce the PBW spanning set of Feigin-Stoyanovsky's subspaces of level $1$ and level $2$ modules to a basis. Monomials from this basis satisfy certain combinatorial conditions, i.e., the difference and initial conditions.

\begin{prop}
	Monomials that do not satisfy difference conditions for level $k$ can be excluded from the spanning set \eqref{E: PBW base def}. 
\end{prop}

\begin{dokaz}
	We show that any monomial $\mm$ that does not satisfy difference conditions necessarily contains a submonomial $\mm'$ of length $k+1$ which is the minimal monomial appearing in some relation. We refer to such minimal monomials as \textit{leading terms} of respective relations. Hence $\mm'$ can be replaced by some greater monomials and, by \eqref{OrdMult_rel} so can $\mm$ also be replaced.
	We show this first for conditions \ref{DC2} and \ref{DC5}, and then give a sketch of the argument for the rest of the cases.
	
	\smallskip
	\ref{DC2} \quad	
	Assume that $\mm$ doesn't satisfy difference condition \ref{DC2}. This means that $\mm$ has a submonomial 
	$$\mm'=x_2(-j-1)^{b_2} \xtu(-j)^{a_\triu} \xcu(-j)^{a_\cetu} x_4(-j)^{a_4} x_3(-j)^{a_3} x_2(-j)^{a_2},$$
	where $b_2+ a_\triu+ a_\cetu+ a_4+ a_3+ a_2 = k+1$.
	To obtain a relation for $\mm'$, start from a relation
	\begin{equation} \label{DC2_rel}
		x_4(z)^{a_\triu +a_4} x_3(z)^{a_\cetu} x_2(z)^{b_2+a_3+a_2}=0
	\end{equation}
	(see \eqref{VOAdc_rel}). 
	
	Denote by $E=x_{\epsilon_3-\epsilon_2}$, $A=x_{\epsilon_4-\epsilon_3}$ and $B=x_{-\epsilon_4-\epsilon_3}$ Chevalley generators of $\mathfrak{l}_0$. The adjoint action of $\mathfrak{l}_0$ on $\g_1$ (and hence on $\gk_1$) is given by the following graph: 
	
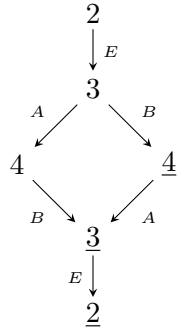
\begin{center}
		\begin{tikzpicture}[node distance=1cm,on grid,>=stealth,bend angle=45,auto,
				every place/.style= {minimum size=6mm,thick,draw=blue!75,fill=blue!20},
				every transition/.style={draw=none,fill=none},
				red place/.style= {place,draw=red!75,fill=red!20},
				every label/.style= {red}]
			
			\node (0) {};
			\node [transition] (3) [above=of 0] {$
					3$};
			\node [transition] (2) [above=of 3] {$
					2$}
			edge [post] node {$\scriptscriptstyle E$} (3) ;
			\node [transition] (4) [left=of 0] {$
					4$}
			edge [pre] node {$\scriptscriptstyle A$} (3);
			\node [transition] (44) [right=of 0] {$
					\underline{4}$}
			edge [pre] node[swap] {$\scriptscriptstyle B$} (3);
			\node [transition] (33) [below=of 0] {$
					\underline{3}$}
			edge [pre] node[swap] {$\scriptscriptstyle A$} (44)
			edge [pre] node {$\scriptscriptstyle B$} (4);
			\node [transition] (22) [below=of 33] {$
					\underline{2}$}
			edge [pre] node {$\scriptscriptstyle E$} (33);
		\end{tikzpicture}
	\end{center}
	With a proper normalization, this corresponds to the following relations:
	\begin{eqnarray*}
		& & \hspace{-2ex} E\cdot x_2=x_3,\quad E\cdot \xtu=\xdu,\quad E\cdot x_3=E\cdot x_4=E\cdot \xcu=E\cdot \xdu=0, \\
		& & \hspace{-2ex} A\cdot x_3=x_4,\quad A\cdot \xcu=\xtu,\quad A\cdot x_2=A\cdot x_4=A\cdot \xtu=A\cdot \xdu=0, \\
		& & \hspace{-2ex} B\cdot x_3=\xcu,\quad B\cdot x_4=\xtu,\quad B\cdot x_2=B\cdot \xcu =B\cdot \xtu=B\cdot \xdu=0. 
	\end{eqnarray*}
	
	Action on relation \eqref{DC2_rel} by $B^{a_\cetu + a_\triu}$ gives
	\begin{equation} \label{DC2korak2_rel}
		\sum C_{r} \xtu(z)^{a_\cetu+a_\triu-r} \xcu(z)^{r} x_4(z)^{a_4+r-a_\cetu} x_3(z)^{a_\cetu-r} x_2(z)^{b_2+a_3+a_2}=0,
	\end{equation}
	where $C_{r}$ are nonzero integers and the sum goes over $0\leq r \leq a_\cetu$, $a_\cetu-a_4 \leq r$.	
	Next, action on \eqref{DC2korak2_rel} by $E^{a_3}$ gives
	\begin{equation} \label{DC2korak3_rel}
		\sum C_{r,s} \xdu(z)^{s} \xtu(z)^{a_\cetu+a_\triu-r-s} \xcu(z)^{r} x_4(z)^{a_4+r-a_\cetu} x_3(z)^{a_\cetu-r+a_3-s} x_2(z)^{b_2+a_2+s}=0,
	\end{equation}
	where $C_{r,s}$ are nonzero integers and the sum goes over $0\leq r \leq a_\cetu$, $a_\cetu-a_4 \leq r$, $0\leq s\leq a_3$, $s\leq a_\triu+a_\cetu-r$.	
	
	The coefficient of $z^{b_2 j+(a_\triu+ a_\cetu+ a_4+ a_3+ a_2)(j-1)}$ gives a relation between monomials. Notice that all monomials in this relation are of length $k+1$ and have the same total degree. Hence the minimal monomial in the relation will be of {\em the minimal shape}, meaning that the degrees of all factors differ for at most $1$. Equivalently, the minimal monomial will have $b_2$ factors of degree $-j-1$ and $a_\triu+ a_\cetu+ a_4+ a_3+ a_2$ factors of degree $-j$. We need to find the minimal configuration of colors from \eqref{DC2korak3_rel} for this shape, first for the $(-j)$-part of the monomial, and then for the $(-j-1)$-part. Hence the minimal monomial looks like
	$$x_2(-j-1)^{b_2} \xdu(-j)^{s} \xtu(-j)^{a_\cetu+a_\triu-r-s} \xcu(-j)^{r} x_4(-j)^{a_4+r-a_\cetu} x_3(-j)^{a_\cetu-r+a_3-s} x_2(-j)^{a_2+s}.$$
	We look for $r,s$ for which such monomial is minimal. This means that the exponent of $x_2(-j)$ needs to be the least possible, after which the exponent of $x_3(-j)$ needs to be the least possible, and so forth. This first implies that $s$ is minimal, hence $s=0$. Next, $r$ needs to be maximal, hence $r=a_\cetu$. We get that the minimal monomial of the relation is exactly $\mm'$. 
	
	\smallskip
	\ref{DC5} \quad	
	Assume that $\mm$ doesn't satisfy difference condition \ref{DC5}. This means that $\mm$ has a submonomial 
	$$\mm'=x_4(-j-1)^{b_4} x_3(-j-1)^{b_3} x_2(-j-1)^{b_2} \xdu(-j)^{a_\dvau} \xtu(-j)^{a_\triu} x_4(-j)^{a_4},$$
	where $b_4+ b_3+ b_2+ a_\dvau+ a_\triu+ a_4 = k+1$.
	This time start from
	\begin{equation} \label{DC5_rel} x_4(z)^{b_4+a_\dvau +a_\triu + a_4} x_3(z)^{b_3} x_2(z)^{b_2}=0. \end{equation}
	
	Acting on relation \eqref{DC5_rel} by $B^{a_\dvau + a_\triu}$ gives
	\begin{equation} \label{DC5korak2_rel}
		\sum C_{r}\xtu(z)^{a_\triu+a_\dvau - r} \xcu(z)^r x_4(z)^{b_4+a_4+r} x_3(z)^{b_3-r} x_2(z)^{b_2}=0,
	\end{equation}
	where $C_{r}$ are nonzero integers and the sum goes over $0\leq r \leq a_\triu+a_\dvau$, $r\leq b_3$.	
	Next, acting on \eqref{DC5korak2_rel} by $E^{a_\dvau}$ gives
	\begin{equation} \label{DC5korak3_rel}
		\sum C_{s,r} \xdu(z)^{a_\dvau-s} \xtu(z)^{a_\triu- r+s} \xcu(z)^r x_4(z)^{b_4+a_4+r} x_3(z)^{b_3-r+s} x_2(z)^{b_2-s}=0,
	\end{equation}
	where $C_{s,r}$ are nonzero integers and the sum goes over $0\leq r \leq a_\triu+a_\dvau$, $r\leq b_3$, $0\leq s\leq b_2$, $s\leq a_\dvau$.
	
	The coefficient of $z^{(b_4+b_3+b_2 )j+(a_\dvau+a_\triu+ a_4)(j-1)}$ gives a relation between monomials. Again, all monomials in this relation are of length $k+1$ and have the same total degree. Hence the minimal monomial in the relation will be of the minimal shape --it will have $b_4+b_3+b_2$ factors of degree $-j-1$ and $a_\dvau+ a_\triu+ a_4$ factors of degree $-j$. We need to find the minimal configuration of colors from \eqref{DC5korak3_rel} for this shape, first for the $(-j)$-part of the monomial, and then for the $(-j-1)$-part. Hence the minimal monomial looks like
	$$x_4(-j-1)^{b_4+r} x_3(-j-1)^{b_3-r+s} x_2(-j-1)^{b_2-s} \xdu(-j)^{a_\dvau-s} \xtu(-j)^{a_\triu- r+s} \xcu(-j)^r x_4(-j)^{a_4}.$$
	We look for $r,s$ for which such monomial is minimal. Note first that $r$ needs to be minimal, hence $r=0$. Next, $s$ needs to be minimal, hence $s=0$. This gives that the minimal monomial of the relation is exactly $\mm'$. 
	
	For other difference conditions we give only a sketch of the argument.
	
	\smallskip 
	\ref{DC1} \quad	
	Let $b_2+ a_\dvau+ a_\triu+ a_\cetu+ a_4+ a_3 = k +1$. Action of $B^{a_\cetu}E^{a_\dvau+a_\cetu+a_3}$ on the relation
	$$ \xtu(z)^{a_\dvau+a_\triu} x_4(z)^{a_4} x_2(z)^{b_2+a_\cetu+a_3} =0 $$
	gives a relation
	$$\sum C_{s,r} \xdu(z)^{a_\dvau+a_3+a_\cetu-r} \xtu(z)^{a_\triu -a_3-a_\cetu + r+s} \xcu(z)^{a_\cetu-s} x_4(z)^{a_4-s} x_3(z)^{r-a_\cetu+s} x_2(z)^{b_2+a_3+a_\cetu-r}=0,$$
	where $C_{s,r}$ are nonzero integers and the sum goes over 
	$0\leq r\leq a_\dvau+a_3+a_\cetu $, $ a_\triu -a_3 -a_\cetu \leq r \leq b_2 +a_3 +a_\cetu$, 
	$0\leq s \leq a_4$, $a_\cetu -r \leq s\leq a_\cetu$.
	
	\smallskip 
	\ref{DC3} \quad	
	Let $b_3+ b_2+ a_\dvau+ a_\triu+ a_\cetu+ a_4 = k +1$. Action of $E^{a_\dvau}B^{a_\dvau+a_\triu+a_\cetu} $ on the relation
	$$x_4(z)^{a_4+a_\dvau+a_\triu} x_3(z)^{b_3+a_\cetu} x_2(z)^{b_2} =0 $$
	gives a relation
	$$\sum C_{s,r} \xdu(z)^{a_\dvau-s} \xtu(z)^{a_\triu + s + a_\cetu - r} \xcu(z)^r x_4(z)^{a_4+r-a_\cetu} x_3(z)^{b_3+a_\cetu-r+s} x_2(z)^{b_2-s}=0,$$
	where $C_{s,r}$ are nonzero integers and the sum goes over 
	$0\leq r\leq a_\dvau+a_\triu+a_\cetu $, $ r \leq b_3 +a_4$, 
	$0\leq s \leq a_\dvau$.
	
	\smallskip 
	\ref{DC4} \quad	
	Let $b_3+ b_2+ a_\dvau+ a_\cetu+ a_4+ a_3 = k +1$. Action of $A^{a_4}E^{b_3+a_\dvau+a_4+a_3} $ on the relation
	$$\xtu(z)^{a_\dvau} \xcu(z)^{a_\cetu} x_2(z)^{b_3+b_2+a_3+a_4} =0 $$
	gives a relation
	$$\sum C_{s,r} \xdu(z)^{a_\dvau-r} \xtu(z)^{ s + r} \xcu(z)^{a_\cetu -s} x_4(z)^{a_4 - s} x_3(z)^{b_3 + a_3 +r+s} x_2(z)^{b_2-r}=0,$$
	where $C_{s,r}$ are nonzero integers and the sum goes over 
	$0\leq r\leq a_\dvau$, 
	$0\leq s \leq a_4$.
	
	\smallskip 
	\ref{DC6} \quad	
	Let $b_\cetu+ b_3+ b_2+ a_\dvau+ a_\triu+ a_\cetu = k +1$. Action of $E^{a_\dvau}A^{a_\dvau+a_\triu} $ on the relation
	$$\xcu(z)^{b_\cetu+a_\dvau+a_\triu+a_\cetu} x_3(z)^{b_3} x_2(z)^{b_2} =0 $$
	gives a relation
	$$\sum C_{s,r} \xdu(z)^{a_\dvau-s} \xtu(z)^{a_\triu- r+s} \xcu(z)^{b_\cetu+a_\cetu+r} x_4(z)^{r} x_3(z)^{b_3-r+s} x_2(z)^{b_2-s}=0,$$
	where $C_{s,r}$ are nonzero integers and the sum goes over 
	$0\leq r\leq a_\dvau+a_\triu$, 
	$0\leq s \leq a_\dvau$, $r -a_\triu \leq s\leq b_2$.
	
	\smallskip 
	\ref{DC7} \quad	
	Let $b_\cetu+ b_4+ b_3+ b_2+ a_\dvau+ a_\triu = k +1$. Action of $E^{a_\dvau}B^{b_\cetu+a_\dvau+a_\triu} $ on the relation
	$$x_4(z)^{b_4+a_\dvau+a_\triu} x_3(z)^{b_3+b_\cetu} x_2(z)^{b_2} =0 $$
	gives a relation
	$$\sum C_{s,r} \xdu(z)^{a_\dvau-s} \xtu(z)^{b_\cetu+a_\triu- r+s} \xcu(z)^{r} x_4(z)^{b_4-b_\cetu+r} x_3(z)^{b_3+b_\cetu-r+s} x_2(z)^{b_2-s}=0,$$
	where $C_{s,r}$ are nonzero integers and the sum goes over 
	$0\leq r\leq b_\cetu+a_\dvau+a_\triu$, $b_{\cetu}-b_4 \leq r \leq b_3+b_\cetu$, 
	$0\leq s \leq a_\dvau$, $r- b_\cetu -a_\triu \leq s\leq b_2$.
	
	\smallskip 
	\ref{DC8} \quad	
	Let $b_\triu+ b_\cetu+ b_4+ b_2+ a_\dvau+ a_\triu = k +1$. Action of $ E^{a_\dvau}B^{b_\cetu+b_\triu+a_\dvau+a_\triu}$ on the relation
	$$ x_4(z)^{b_\triu+b_4+a_\dvau+a_\triu} x_3(z)^{b_\cetu} x_2(z)^{b_2} =0$$
	gives a relation
	$$\sum C_{s,r} \xdu(z)^{a_\dvau-s} \xtu(z)^{b_\triu+b_\cetu+a_\triu- r+s} \xcu(z)^{r} x_4(z)^{b_4-b_\cetu+r} x_3(z)^{b_\cetu-r+s} x_2(z)^{b_2-s}=0,$$
	where $C_{s,r}$ are nonzero integers and the sum goes over 
	$0\leq r\leq b_\cetu$, $b_{\cetu}-b_4 \leq r$, 
	$0\leq s \leq a_\dvau$, $r- b_\cetu -b_\triu-a_\triu \leq s\leq b_2$.
	
	\smallskip 
	\ref{DC9} \quad	
	Let $b_\triu+ b_\cetu+ b_4+ b_3+ b_2+ a_\dvau = k +1$. Action of $B^{b_\cetu}E^{b_\cetu+b_3+a_\dvau} $ on the relation
	$$x_4(z)^{b_4} x_\triu(z)^{b_\triu+a_\dvau} x_2(z)^{b_\cetu+b_3+b_2} =0 $$
	gives a relation
	$$\sum C_{s,r} \xdu(z)^{b_\cetu+b_3+a_\dvau-s} \xtu(z)^{b_\triu-b_3-b_\cetu + r+s} \xcu(z)^{b_\cetu-r} x_4(z)^{b_4-r} x_3(z)^{r+s - b_\cetu} x_2(z)^{b_2 +b_3 +b_\cetu -s}=0,$$
	where $C_{s,r}$ are nonzero integers and the sum goes over 
	$0\leq s \leq b_3 + b_\cetu+a_\dvau$, $b_3 + b_\cetu -b_\triu \leq s\leq b_2 +b_3 +b_\cetu$, 
	$0\leq r\leq b_\cetu$, $b_{\cetu}-s \leq r \leq b_4$.
	
	\smallskip 
	\ref{DC10} \quad	
	Let $b_\dvau+ b_\triu+ b_\cetu+ b_4+ b_3+ a_\dvau = k+1$. Action of $B^{b_\cetu} E^{b_\dvau+b_\cetu+b_3+a_\dvau} $ on the relation
	$$x_\triu(z)^{b_\dvau+b_\triu+a_\dvau} x_4(z)^{b_4} x_2(z)^{b_\cetu+b_3} =0 $$
	gives a relation
	$$\sum C_{s,r} \xdu(z)^{b_\dvau+b_\cetu+b_3+a_\dvau-s} \xtu(z)^{b_\triu-b_3-b_\cetu + r+s} \xcu(z)^{b_\cetu-r} x_4(z)^{b_4-r} x_3(z)^{r+s - b_\cetu} x_2(z)^{b_3 +b_\cetu -s}=0,$$
	where $C_{s,r}$ are nonzero integers and the sum goes over 
	$0\leq s \leq b_3 + b_\cetu$, $s\geq b_3+b_\cetu-b_\triu$, 
	$0\leq r\leq b_\cetu$, $b_{\cetu}-s\leq r \leq b_4$.
\end{dokaz}

We also define initial conditions for a monomial: let $b_t$ (resp. $b_{\underline{t}}$) denote the number of $x_t(-1)$ (resp. $x_{\underline{t}}(-1)$) factors in $\mm$. We say that the monomial $\mm$ satisfies \emph{initial conditions} for $L(\Lambda)$, $\Lambda=k_0\Lambda_0+ k_1 \Lambda_1+k_2 \Lambda_2 + k_3 \Lambda_3 + k_4 \Lambda_4$ if the following frequency conditions hold:

\begin{enumerate}[label=(\roman*)]
	\item \label{IC1} $b_2 \leq k_0$,
	\item \label{IC2} $b_3+ b_2 \leq k_0+k_2$,
	\item \label{IC3} $b_4+ b_3+ b_2 \leq k_0+k_2+k_3$,
	\item \label{IC4}$b_\cetu+ b_3+ b_2\leq k_0+k_2+k_4$,
	\item \label{IC5}$b_\cetu+ b_4+ b_3+ b_2 \leq k_0+k_2+k_3+k_4$,
	\item \label{IC6}$b_\triu+ b_\cetu+ b_4+ b_2 \leq k_0+k_2+k_3+k_4$,
	\item \label{IC7}$b_\triu+ b_\cetu+ b_4+ b_3+ b_2\leq k_0+2k_2+k_3+k_4$,
	\item \label{IC8}$b_\dvau+ b_\triu+ b_\cetu+ b_4+ b_3 \leq k_0+2k_2+k_3+k_4$.
\end{enumerate}

\begin{remark} \label{IC0_remark} {\em
		As in other cases (cf. \cite{Ba, BPT, T2}), initial conditions can be expressed in terms of difference
		conditions by adding ``imaginary'' $(0)$-factors to $\mm$. Let
		$$\mm_\Lambda=\xdu(0)^{k_1} \xtu(0)^{k_2} \xcu(0)^{k_3} x_4(0)^{k_4} x_3(0)^{k_2} x_2(0)^{k_1}.$$
		Then a monomial $\mm$ satisfies difference and initial conditions for $L(\Lambda)$ if
		and only if $\mm'=\mm \cdot \mm_\Lambda$ satisfies difference conditions. In
		fact, initial conditions are defined in this way so that the property holds.
	}\end{remark}

\begin{prop}
	Monomials that do not satisfy initial conditions for level $k$ can be excluded from the spanning set \eqref{E: PBW base def}. 
\end{prop}

\begin{dokaz}
	We show that for a monomial 
	$\mm$ that doesn't satisfy initial conditions either $\mm v_\Lambda=0$, or the monomial vector $\mm v_\Lambda$ is minimal for some relation. Hence $\mm v_\Lambda$ can be replaced by greater monomial vectors.
	
	Assume that $\mm=x_4 (-1)^{b_4} x_3 (-1)^{b_3} x_2 (-1)^{b_2}$ doesn't satisfy initial condition \ref{IC3}, i.e. that $b_4+b_3+b_2>k_0+k_2+k_3$. Then, by Lemma \ref{IC1_lm}, 
	$$\mm v_\Lambda = \mm \left(v_{\Lambda_{0}}^{\otimes k_0} \otimes v_{\Lambda_{1}}^{\otimes k_1} \otimes v_{\Lambda_{2}}^{\otimes k_2} \otimes v_{\Lambda_{3}}^{\otimes k_3} \otimes v_{\Lambda_{4}}^{\otimes k_4}\right) = 0,$$
	because every factor of $\mm$ annihilates $v_{\Lambda_1}$. So 
	all monomial vectors that contain $\mm$ as a submonomial can be excluded from the spanning set. Similar argument can be given if $\mm$ doesn't satisfy initial condition \ref{IC1}, \ref{IC2} and \ref{IC4}.
	
	Next, assume that $\mm= \xtu(-1)^{b_\triu} \xcu(-1)^{b_\cetu} x_4(-1)^{b_4} x_3(-1)^{b_3} x_2(-1)^{b_2} $ doesn't satisfy initial condition \ref{IC7}. 
	Set $\Lambda'=k_0 \Lambda_0+k_2 \Lambda_2 +k_3 \Lambda_3+k_4 \Lambda_4$ and $\Lambda''=k_1 \Lambda_1$. Then, by Lemma \ref{IC1_lm},
	$$\mm v_\Lambda =\mm \left(v_{\Lambda'}\otimes v_{\Lambda''}\right) =\left(\mm v_{\Lambda'}\right)\otimes v_{\Lambda''}.$$
	The claim now follows from the relation for difference condition \ref{DC9} for the module $L(\Lambda')$.
	Similar argument can be given if $\mm$ doesn't satisfy initial condition \ref{IC8}.
	
	Finally, assume that $\mm= \xcu(-1)^{b_\cetu} x_4(-1)^{b_4} x_3(-1)^{b_3} x_2(-1)^{b_2} $ doesn't satisfy initial condition \ref{IC5}.
	Action of $B^{b_\cetu}$ on the relation
	$$x_4(-1)^{b_4} x_3(-1)^{b_\cetu+b_3} x_2(-1)^{b_2} v_\Lambda=0$$
	gives the following relation
	$$\sum C_{r} \xtu(-1)^{r} \xcu(-1)^{b_\cetu-r} x_4(-1)^{b_4-r} x_3(-1)^{b_3 + r} x_2(-1)^{b_2} v_\Lambda=0,$$
	where $C_{r}$ are nonzero integers and the sum goes over 
	$0\leq r\leq b_\cetu$, $ r \leq b_4$.
	The monomial vector $\mm v_\Lambda$ is the minimal one in this relation, so it can be excluded from the spanning set.
	Similar argument can be given if $\mm$ doesn't satisfy initial condition \ref{IC6}.
\end{dokaz}

\begin{remark} {\em
		The method used in \cite{BPT, T2} for proving linear independence of the spanning set for $W(\Lambda)$ in the level $k$ case relies crucially on the fact that every monomial satisfying IC and DC on level $k$ admits a factorization with respect to the tensor factors that satisfies respective IC and DC. In the case of $D_4$ this doesn't hold: a monomial
		
		\[
			x_3(-1)\,x_{\underline{3}}(-1)\;
			x_2(-2)\,x_3(-2)\;
			x_2(-3)\,x_3(-3)\,x_{\underline{3}}(-3)\;
			x_3(-4)^2
		\]
		satisfies DC for level 3 modules, but it does not admit a factorization on two monomials satisfying DC for level 1 and level 2 modules, respectively.
		Thus the logic of the linear independence proof from \cite{BPT, T2} cannot be applied here.
	}\end{remark}

\begin{remark}{\em
		We compute the number of leading terms corresponding to the following critical equalities:
		\begin{enumerate}[label=(\roman*)]
			\item \label{CC1} $b_2+ a_\dvau+ a_\triu+ a_\cetu+ a_4+ a_3 = k + 1$,
			\item \label{CC2} $b_2+ a_\triu+ a_\cetu+ a_4+ a_3+ a_2 = k + 1$,
			\item \label{CC3} $b_3+ b_2+ a_\dvau+ a_\triu+ a_\cetu+ a_4 = k + 1$,
			\item \label{CC4} $b_3+ b_2+ a_\dvau+ a_\cetu+ a_4+ a_3 = k + 1$,
			\item \label{CC5} $b_4+ b_3+ b_2+ a_\dvau+ a_\triu+ a_4 = k + 1$,
			\item \label{CC6} $b_\cetu+ b_3+ b_2+ a_\dvau+ a_\triu+ a_\cetu = k + 1$,
			\item \label{CC7} $b_\cetu+ b_4+ b_3+ b_2+ a_\dvau+ a_\triu = k + 1$,
			\item \label{CC8} $b_\triu+ b_\cetu+ b_4+ b_2+ a_\dvau+ a_\triu = k + 1$,
			\item \label{CC9} $b_\triu+ b_\cetu+ b_4+ b_3+ b_2+ a_\dvau = k + 1$,
			\item \label{CC10} $b_\dvau+ b_\triu+ b_\cetu+ b_4+ b_3+ a_\dvau = k + 1$.
		\end{enumerate}
		Note that these leading terms are precisely the minimal monomials, with respect to order defined in Section \ref{S: FS}, that violate initial and difference conditions.
		
		Let $\#LT(k)$ denote the number of leading terms at level $k$. In order to obtain $\#LT(k)$, we count the number of different solutions for the critical equalities \cref{CC1,CC2,CC3,CC4,CC5,CC6,CC7,CC8,CC9,CC10}.
		Equation \ref{CC1} gives \(\binom{k+6}{5}\) solutions.
		A solution of \ref{CC2} differs from the previous ones if $a_2 > 0$; this gives another  $\binom{k+5}{5}$ solutions.
		A solution of \ref{CC3} is new only  if $b_3 > 0$; this gives another \(\binom{k+5}{5}\) solutions.
		Proceed analogously as above for the rest of the equations. One gets
		\begin{equation*} 
			\mathrm{\#}\mathrm{LT}(k) = \binom{k+6}{5} + 6\binom{k+5}{5} + 3\binom{k+4}{5} = \frac{(k+2)^2(k+3)^2(k+4)}{12}
		\end{equation*} as the total number of leading terms at level \(k\).
		
		We can compare $\#LT(k)$ with dimension of the $\mathfrak{l}_0$-module $V_k$ spanned by $x_\theta(-1)^{k+1}$. By Weyl dimension formula (cf. \cite{H}) 
		\begin{equation*}
			\dim V_k =\frac{(k+2)(k+3)^2(k+4)}{12}.
		\end{equation*}
		Hence
		\[
			\mathrm{\#}\mathrm{LT}(k) = (k+2) \cdot \dim V_k.
		\]
		
		For every relation there are $k+2$ leading terms, each corresponding to different minimal shapes. Hence, the relations \cref{CC1,CC2,CC3,CC4,CC5,CC6,CC7,CC8,CC9,CC10} describe all the leading terms corresponding to the $\mathfrak{l}_0$ action on the Frenkel-Kac-Segal relation
		\[
			x_{\theta}(z)^{k+1} = 0.
		\]
		This gives a reason to believe that the spanning set described above is in fact a basis.	
	}\end{remark}

\section{Spanning sets of standard modules} \label{S: Spanning}

As in other cases (cf. \cite{FS, P1, P2, P3}), a spanning set of the whole standard module $L(\Lambda)$ can be obtained as a limit of a spanning set for $W(\Lambda)$.

Set
$$e=\prod_{\gamma\in\Gamma} e_\gamma.$$

\begin{prop}[\cite{P1, P2, P3}] \label{StModGenSkup_prop}
	$L(\Lambda)=\langle e \rangle U(\gt_1) v_\Lambda$.
\end{prop}

Proposition \ref{StModGenSkup_prop} was proven in the case of $A_\ell$ and $B_2$ in \cite{P1, P2, P3}, and an analogous proof holds in this case.

Since
$$		\gamma_2 + \gamma_3+ \gamma_4 + \gamma_\cetu+\gamma_\triu + \gamma_\dvau = 6\epsilon_1= 6 \omega,
$$
then $e= C [\omega]^{6}$, for some $C\neq 0$. Hence

\begin{kor} \label{omegaWgen}
	$L(\Lambda)=\langle [\omega]^2 \rangle U(\gt_1) v_\Lambda$.
\end{kor}

Define extremal vectors
$$v_{\Lambda}^{(-2m)}=[\omega]^{-2m} v_{\Lambda},\qquad m\in\Z$$
and the corresponding shifted Feigin-Stoyanovsky's type subspaces
$$W_{-2m}= U(\gt_1)v_{\Lambda}^{(-2m)}.$$
This gives a sequence of embeddings
\begin{eqnarray} 
	\label{extremal vec}
	&	\cdots \xrightarrow{[\omega]^{-2}} W_0 \xrightarrow{[\omega]^{-2}} W_{-2} \xrightarrow{[\omega]^{-2}} W_{-4}\xrightarrow{[\omega]^{-2}} W_{-6} \xrightarrow{[\omega]^{-2}}\ \cdots & \\
	\nonumber
	&	\cdots \xrightarrow{[\omega]^{-2}} v_{\Lambda} \xrightarrow{[\omega]^{-2}} v_{\Lambda}^{(-2)} \xrightarrow{[\omega]^{-2}} v_{\Lambda}^{(-4)}\xrightarrow{[\omega]^{-2}} v_{\Lambda}^{(-6)} \xrightarrow{[\omega]^{-2}}\ \cdots & 
\end{eqnarray}

By Corollary \ref{omegaWgen}, 
$$L(\Lambda)=\bigcup_{m\in\Z} W_{-2m}.$$

Since
$$[\omega]^{-2m} {\bf m} v_{\Lambda}={\bf m}^{+2m} v_{\Lambda}^{(-2m)},$$
a spanning set of $W_{-2m}$ can be obtained by taking a spanning set of $W$ and shifting degrees of monomials. We say that
${\bf m}$ satisfies {\em initial conditions} $IC_{2m}$ for $W(\Lambda)$ if ${\bf m}^{+2m}$ satisfies initial conditions, and that is if and only if it doesn't contain elements of degree $2m$ or greater  and elements of degree $2m-1$ satisfy additional conditions corresponding to initial conditions for elements of degree $-1$ in ${\bf m}^{+2m}$.

\begin{prop} The set
	$$\mathcal{B}_{-2m}=\{{\bf m}v_{\Lambda}^{(-2m)} \,|\, {\bf m} \ \textrm{satisfies}\ DC\ \textrm{and}\ IC_{2m}\ \textsf{ for}\ W(\Lambda)\},$$
	is a generating set for $W_{-2m}$.
\end{prop}

To describe the sequence of embeddings \eqref{extremal vec}, set 

$${\bf n}_0= x_2(-1) x_\dvau(-1), \qquad {\bf n}_1= 1, \qquad {\bf n}_2= x_3(-1) x_\triu (-1),$$
$$ {\bf n}_3= x_4(-1),  \qquad{\bf n}_4= x_\cetu(-1). $$
Furthermore, set
$${\bf m}_0= {\bf n}_1^- {\bf n}_0, \qquad{\bf m}_1= {\bf n}_0^- {\bf n}_1, \qquad{\bf m}_2= {\bf n}_2^- {\bf n}_2,$$
$${\bf m}_3= {\bf n}_4^- {\bf n}_3, \qquad{\bf m}_4= {\bf n}_3^- {\bf n}_4. $$

Lemmas \ref{DC1_lm} and \ref{IC1_lm}, relations \cref{SCrel1,SCrel2,SCrel3,SCrel4,SCrel5} and Remark \ref{IC0_remark} give

\begin{prop} \label{maxarg}
	\begin{enumerate}[label=(\roman*)]
		\item Monomial ${\bf n}_i$, $i=0,1,2,3,4$, is the maximal monomial with factors of degree greater than $-1$ that acts nontrivially on $v_{\Lambda_i}$ and
		      $$ {\bf n}_0 v_{\Lambda_0}=C_0 [\omega] v_{\Lambda_{1}},\qquad {\bf n}_1 v_{\Lambda_1}=C_1 [\omega] v_{\Lambda_{0}},\qquad {\bf n}_2 v_{\Lambda_2}=C_2 [\omega] v_{\Lambda_{2}},$$ $${\bf n}_3 v_{\Lambda_3}=C_3 [\omega] v_{\Lambda_{4}},\qquad {\bf n}_4 v_{\Lambda_4}=C_4 [\omega] v_{\Lambda_{3}},$$ for some $C_i\neq 0, i=0,\dots,4$.
		\item Monomial ${\bf m}_i$, $i=0,1,2,3,4$, is the maximal monomial with factors of degree greater than $-2$ that acts nontrivially on $v_{\Lambda_i}$
		      and
		      $${\bf m}_iv_{\Lambda_i}=C_i'[\omega]^2 v_{\Lambda_i},$$
		      for some $C_i'\neq 0, i=0,\dots,4$.
		\item A monomial ${\bf m}$ satisfies difference and initial conditions for $W(\Lambda_i)$ if and only if ${\bf m} {\bf m}_i^{+2}$ satisfies DC.
	\end{enumerate}
\end{prop}

For a higher level module $L(\Lambda)$, $\Lambda=k_0 \Lambda_0 + k_1 \Lambda_1 + k_2 \Lambda_2 + k_3 \Lambda_3 + k_4 \Lambda_4$, with the highest weight vector $v_\Lambda=v_{\Lambda_{0}}^{\otimes k_0} \otimes v_{\Lambda_{1}}^{\otimes k_1} \otimes v_{\Lambda_{2}}^{\otimes k_2} \otimes v_{\Lambda_{3}}^{\otimes k_3} \otimes  v_{\Lambda_{4}}^{\otimes k_4}$ set
$$\mm_\Lambda=\mm_0^{k_0} \mm_1^{k_1} \mm_2^{k_2} \mm_3^{k_3} \mm_4^{k_4}. $$
Then Proposition \ref{maxarg} and \eqref{OrdMult_rel} give
\begin{eqnarray*}
	\mm_\Lambda v_\Lambda & = & (\mm_0 v_{\Lambda_{0}})^{\otimes k_0} \otimes (\mm_1 v_{\Lambda_{1}})^{\otimes k_1} \otimes (\mm_2 v_{\Lambda_{2}})^{\otimes k_2} \otimes (\mm_3 v_{\Lambda_{3}})^{\otimes k_3} \otimes (\mm_4 v_{\Lambda_{4}})^{\otimes k_4} \\
	& = & (C_0' [\omega]^2 v_{\Lambda_{0}})^{\otimes k_0} \otimes (C_1' [\omega]^2 v_{\Lambda_{1}})^{\otimes k_1} \otimes (C_2' [\omega]^2 v_{\Lambda_{2}})^{\otimes k_2} \otimes (C_3' [\omega]^2 v_{\Lambda_{3}})^{\otimes k_3} \otimes (C_4' [\omega]^2 v_{\Lambda_{4}})^{\otimes k_4} \\
	& = & C' [\omega]^2 v_\Lambda, 
\end{eqnarray*}
where $ C', C_1',\dots,C_4'\neq 0$. 

Note that, up to a scalar, 
$$\mm v_{\Lambda}	= \mm [\omega]^{-2} [\omega]^2 v_{\Lambda} = \mm \mm_{\Lambda}^{+2} v_{\Lambda}^{(-2)}, $$
for a monomial $\mm$. So, the spanning set of $W$ can be embedded into the spanning set of $W_{-2}$, 
$$\mathcal{B}_0 \subset \mathcal{B}_{-2}.$$
Proceed inductively to obtain a chain of inclusions
\begin{equation}\label{B_mbeddings}
	\mathcal{B}_0 \subset \mathcal{B}_{-2}\subset \mathcal{B}_{-4}\subset \cdots \subset \mathcal{B}_{-2t}\subset \cdots.
\end{equation}

To obtain a spanning set of the whole $L(\Lambda)$ one should take ``the limit $t\to\infty$'', i.e. take inductive limit of the sequence \eqref{B_mbeddings}. Formally, let 
$$v_\Lambda^{(-2t)} = \mm_\Lambda^{+2t+2} \mm_{\Lambda}^{+2t+4}\cdots v_{\Lambda}^{(-\infty)}.$$
A semi-infinite monomial $\mm$ is said to have a periodic tail, or that it {\em stabilizes}, if from some point on it consists of successive shifts of $\mm_\Lambda$,
$$\mm=\mm' \mm_{\Lambda}^{+2t+2} \mm_{\Lambda}^{+2t+4}\cdots,$$
for some $t\in\Z$.

\begin{tm} The set
	\begin{equation*}
		\overline{\mathcal{B}} = 
		\{ \mm v_{\Lambda}^{(-\infty)} \,|\, 
		\mm \ \textrm{stabilizes and satisfies}\ DC\}
	\end{equation*}
	is a generating set for $L(\Lambda)$.
\end{tm}

From \cite{Ba} we know that $\mathcal{B}_0$ is a basis of $W(\Lambda)$ for $\Lambda$ a fundamental weight, hence

\begin{kor} If $\Lambda$ is a fundamental weight, the set
	\begin{equation*}
		\overline{\mathcal{B}} = 
		\{ \mm v_{\Lambda}^{(-\infty)} \,|\, 
		\mm \ \textrm{stabilizes and satisfies}\ DC\}
	\end{equation*}
	is a basis for $L(\Lambda)$.
\end{kor}

\section*{Acknowledgements}

GT was funded by the European Union – NextGenerationEU, under the National Recovery and Resilience Plan, through the project “Fundamental Research in Mathematics” (UFZG-NPOO-25-589).


\end{document}